\documentclass[11pt]{article}
\usepackage{amsfonts,latexsym}
\renewcommand{\theequation}{\arabic{section}.\arabic{equation}}
\newcommand{\cleqn}{\setcounter{equation}{0}}
\newcommand{\clth}{\setcounter{theorem}{0}}
\newcommand {\sectionnew}[1]{\section{#1}\cleqn\clth}
\newcommand{\beq}{\begin{equation}}
\newcommand{\eeq}{\end{equation}}
\newcommand{\beqa}{\begin{eqnarray}}
\newcommand{\eeqa}{\end{eqnarray}}
\newcommand{\beaa}{\begin{eqnarray*}}
\newcommand{\eaa}{\end{eqnarray*}}
\newcommand{\nn}{\hfill\nonumber}
\newcommand{\text}{\textrm}
\newcommand \nc {\newcommand}
\nc \proof {{\em{Proof.\/}} }
\nc \qed {$\Box$\hfill}
\newtheorem{theorem}{Theorem}[section]
\newtheorem{lemma}[theorem]{Lemma}
\newtheorem{proposition}[theorem]{Proposition}
\newtheorem{corollary}[theorem]{Corollary}
\newtheorem{definition}[theorem]{Definition}
\newtheorem{example}[theorem]{Example}
\newtheorem{remark}[theorem]{Remark}
\nc \bth[1] { \begin{theorem}\label{t#1} }
\nc \ble[1] { \begin{lemma}\label{l#1} }
\nc \bpr[1] { \begin{proposition}\label{p#1} }
\nc \bco[1] { \begin{corollary}\label{c#1} }
\nc \bde[1] { \begin{definition}\label{d#1}\rm }
\nc \bex[1] { \begin{example}\label{e#1}\rm }
\nc \bre[1] { \begin{remark}\label{r#1}\rm }
\nc \bcon[1] { \medskip\noindent{\it{Conjecture #1}} }
\nc \bqu[1]  { \medskip\noindent{\it{Question #1}} }
\nc {\eth} { \end{theorem} }
\nc {\ele} { \end{lemma} }
\nc {\epr} { \end{proposition} }
\nc {\eco} { \end{corollary} }
\nc {\ede} { \end{definition} }
\nc {\eex} { \end{example} }
\nc {\ere} { \end{remark} }
\nc {\econ} {\smallskip}
\nc {\equ} {\smallskip}
\nc \eqref[1] {{\rm{(\ref{#1})}}}
\nc \thref[1]{Theorem \ref{t#1}}
\nc \leref[1]{Lemma \ref{l#1}}
\nc \prref[1]{Proposition \ref{p#1}}
\nc \coref[1]{Corollary \ref{c#1}}
\nc \deref[1]{Definition \ref{d#1}}
\nc \exref[1]{Example \ref{e#1}}
\nc \reref[1]{Remark \ref{r#1}}
\newcommand {\normprod}[1]{ {\textrm{:}}{#1}{\textrm{:}} } 
\def \W {W_{1+\infty}}
\def \WN {\W(N)}
\def \a {\alpha}

\def \D {{\mathcal D}}

\def \M {{\mathcal M}}

\def \Cset {{\mathbb C}}
\def \Zset {{\mathbb Z}}
\def \Nset {{\mathbb N}}
\def \Vset {{\mathbb V}}
\def \ord { {\mathrm{ord}} }

\def \span { {\mathrm{span}} }

\def \card { {\mathrm{card}} }
\def \mult { {\mathrm{mult}} }

\def \diag { {\mathrm{diag}} }
\def \res { {\mathrm{Res}} }
\renewcommand \ker { {\mathrm{Ker}} }
\nc \Gr {Gr}
\nc \GRN { \Gr^{(N)} }
\nc \GRA[1] { \Gr_A^{(#1)} }   
\nc \GRAN { \GRA{N} }
\nc \GrA[1] { \Gr_A(#1) }
\nc \GrAa { \GrA{\alpha} }
\nc \GRB[1] { \Gr_B^{(#1)} }   
\nc \GRBN { \GRB{N} }
\nc \GrB[1] { \Gr_B(#1) }
\nc \GrBb { \GrB{\beta} }
\nc \GRMB[1] { \Gr_{MB}^{(#1)} }   
\nc \GRMBN { \GRMB{N} }
\nc \GrMB[1] { \Gr_{MB}(#1) }
\nc \GrMBb { \GrMB{\beta} }

\begin{document}
\title{{\LARGE\bf{
Highest weight modules over
$W_{1+\infty}$ algebra and the bispectral problem
}}}
\author{
B.~Bakalov
\thanks{E-mail: bbakalov@fmi.uni-sofia.bg}
\quad
E.~Horozov
\thanks{E-mail: horozov@fmi.uni-sofia.bg}
\quad
M.~Yakimov
\thanks{E-mail: myakimov@fmi.uni-sofia.bg}
\\ \hfill\\ \normalsize \textit{
Department of Mathematics and Informatics, }\\
\normalsize \textit{Sofia University , 5 J. Bourchier Blvd.,
Sofia 1126, Bulgaria }     }
\date{}
\maketitle
\vspace{-8cm}
\begin{flushright}
{\tt{ q-alg/9602012 }}
\end{flushright}
\vspace{8cm}
\setcounter{section}{-1}
\sectionnew{Introduction}
This paper is the last of a series of papers devoted to the bispectral problem
\cite{BHYc}--\cite{BHY3}. Here we examine the connection between the bispectral
operators constructed in \cite{BHY3} and the Lie algebra $\W$ (and its
subalgebras). To give a more detailed idea of the contents of the present paper
we briefly recall the results of \cite{BHY1}--\cite{BHY3} which we need.

In \cite{BHY1} we built large families of representations of $\W$. For each
$\beta\in\Cset^N$ we defined a tau-function $\tau_\beta(t)$ which we called
{\em{Bessel tau-function\/}}. We proved that it
is a highest weight vector for a representation $\M_\beta$ of the algebra $\W$
with central charge $N$. In \cite{BHY3} we introduced a version of Darboux
transformation, which we called {\em{monomial\/}}, on the corresponding wave
functions $\Psi_\beta(x,z)$  (see also Subsect.~1.2) and showed that the
resulting wave functions are bispectral. For example all bispectral operators
{}from \cite{DG, Z} can be obtained in this way.

The present paper establishes closer connections between $\W$ and the
bispectral problem. Our first result (\thref{7.1}) shows that a tau-function
is a monomial Darboux transformation of a Bessel tau-function if and only if it
belongs to one of the modules $\M_\beta$. This type of connection between the
representation theory (of $\W$) and the bispectral problem is, to the best of
our knowledge, new even for the bispectral tau-functions of
Duistermaat and Gr\"unbaum \cite{DG}.

The second of the questions we try to answer in the present paper originates
{}from Duistermaat and Gr\"unbaum
\cite{DG}. They noticed that their rank 1 bispectral operators are
invariant under the KdV-flows and asked if there is a hierarchy of symmetries
for the rank 2 bispectral operators. The latter question was answered
affirmatively by Magri and Zubelli \cite{MZ} who showed that the algebra
$Vir^+$ (the subalgebra of the Virasoro algebra spanned by the operators of
non-negative weight) is tangent to the manifold of rank 2 bispectral operators.
Here we obtain generalizations of these results as follows.

First we show that the flows generated by $\W^+(N)$ leave our manifold of
monomial Darboux transformations $\GRMBN$ invariant (see \thref{7.7}). An
important feature of our proof is that it naturally follows from the results
contained in \thref{7.1} about the tau-functions in the modules $\M_\beta$.
For this reason
we believe that even in the case $N=2$ \cite{MZ} it gives a better explanation
of the origin of the flows.
Note that the corresponding bispectral operators need not to be of
order $N$ as in \cite{MZ}. Our next result touches upon this particular
situation. We consider the manifolds of rank $N$ polynomial Darboux
transformations (see Subsect.~1.2) of Bessel tau-functions for which the spectral
algebra contains an operator of order $N$. Then a natural bosonic realization
of $Vir^+_N$ generates flows leaving such manifolds invariant (see
\thref{7.8}).
For $N=2$ our theorem coincides with cited above result of \cite{MZ}.

The monomial Darboux transformations form a subfamily of a larger
class of solutions to the bispectral problem -- the polynomial Darboux
transformations of Bessel and Airy planes \cite{BHY3}.
It is a very interesting open problem to find hierarchies of symmetries
preserving the manifolds of polynomial Darboux transformations. We think that
this problem is also connected to representation theory. Perhaps the vertex
operator algebra structure of $\W$ \cite{FKRW} and of (certain completions of)
the modules $\M_\beta$ will help in tackling this question.

{\flushleft{\bf{Acknowledgement}}}

\medskip\noindent
This work was partially supported by Grant MM--523/95 of Bulgarian
Ministry of Education, Science and Technologies.

\sectionnew{Preliminaries}
Here we have collected some facts and notation needed in Sect.\
2 and 3. Most proofs are standard but technical and are given in the Appendix.
\subsection{}
In this subsection we recall some facts and notation from Sato's
theory of KP-hierarchy \cite{S, DJKM, SW} needed in the paper.
We use the approach of V.~Kac and D.~Peterson based on infinite wedge products
(see e.g.\ \cite{KRa}) and the recent survey paper by P.~van~Moerbeke
\cite{vM}.

Consider the infinite-dimensional vector space of formal series
$$
\Vset=\Bigl\{ \sum_{k\in\Zset} a_kv_k \Big|\; a_k=0\ {\rm for}\
k\ll 0\Bigr\}.
$$
Then define the fermionic Fock space $F$ to be the direct sum of the spaces
$F^{(m)}$ (states with a charge number $m$) consisting of formal
infinite sums of semi-infinite wedge monomials
$$
v_{i_0}\wedge v_{i_1}\wedge\ldots
$$
such that $i_0>i_1>\ldots$ and $i_k=m-k$ for $k\gg0$.

There exists a well known linear isomorphism, called a {\em{boson-fermion
correspondence\/}}:
\beq
\sigma\colon F \to
B := \Cset\left[\left[t_1, t_2, \ldots; Q,Q^{-1}\right]\right]
\label{bo-fe}
\eeq
(see \cite{KRa} and the Appendix).

{\em Sato's Grassmannian\/} $Gr$ \cite{S, DJKM, SW} consists of all
subspaces $W\subset\Vset$ which have an admissible basis
$$
w_k=v_k+\sum_{i>k}w_{ik}v_i,\quad k=0,-1,-2,\ldots
$$
To a plane $W\in Gr$ we associate a state
$|W\rangle\in F{}^{(0)}$ as follows
$$
|W\rangle = w_0\wedge w_{-1}\wedge w_{-2}\wedge\ldots
$$
A change of the admissible basis results in a multiplication of $|W\rangle$ by
a non-zero constant. Thus we define an embedding of $Gr$ into the
projectivization of $F{}^{(0)}$ which is called a Pl\"ucker
embedding.
One of the main objects of Sato's theory
is the {\em  tau-function\/} of $W$ defined as the image of
$|W\rangle$ under the boson-fermion correspondence \eqref{bo-fe}
\beq
\tau_W(t) =\sigma(|W\rangle)=\sigma(w_0\wedge w_{-1}\wedge w_{-2}\wedge\ldots).
\label{tau}
\eeq
It is a formal power series in the variables $t_1,t_2,\ldots$, i.e.\
an element of
$B^{(0)} :=
\Cset\left[\left[t_1, t_2, \ldots\right]\right].$
Another important function connected to $W$ is the {\em  Baker\/} or {\em wave
function\/}
\beq
\Psi_W(t,z)= e^{\sum_{k=1}^\infty t_kz^k}
\frac{\tau\left(t-[z^{-1}]\right)}{\tau(t)},
\label{1.7}
\eeq
where $[z^{-1}]$ is the vector $\left(z^{-1}, z^{-2}/2,\ldots\right)$.
Most often $\Psi_W$ is viewed as a formal series. Introducing the vertex
operator \beq
X(t,z)=\exp\left(\sum_{k=1}^\infty t_k z^k\right)
\exp\left(-\sum_{k=1}^\infty \frac1{kz^k}\frac{\partial}{\partial t_k}\right)
\label{1.8}
\eeq
the above formula \eqref{1.7} can be written as
\beq
\Psi_W(t,z)=\frac{X(t,z)\tau(t)}{\tau(t)}.
\label{1.9}
\eeq
We often use the formal series $\Psi_W(x,z) = \Psi_W(t,z)|_{t_1=x,
t_2=t_3=\cdots=0}$, which we call again a Baker function.
The Baker function $\Psi(x,z)$ contains the whole information about $W$ and
hence about $\tau_W$, as the vectors $w_{-k}=\partial^k_x\Psi_W(x,z)|_{x=0}$
form an admissible basis of $W$ (if we take $v_k=z^{-k}$ as a basis of
$\Vset$).

We also use the standard notation
$$
\Gr^{(N)} = \{V\in\Gr | z^N V \subset V\}.
$$
For $V\in\Gr^{(N)}$ there exists an operator $L_V(x,\partial_x)$
of order $N$ such that
$$
L_V(x,\partial_x) \Psi_V(x,z) = z^N \Psi_V(x,z)
$$
and the corresponding tau-function $\tau_V(t)$ does not depend on
$t_N, t_{2N},\ldots$
\subsection{}
Here we shall briefly recall the definition of Bessel wave function and
of monomial Darboux transformations from it. For more details see \cite{BHY3}.

Let $\beta \in \Cset^N$ be such that
\beq
\sum_{i=1}^{N}\beta_i = \frac{N(N-1)}{2}.
\label{2.15}
\eeq
\bde{bess} \cite{F, Z, BHY1}
{\it {Bessel wave function}} is called the unique wave function
$\Psi_\beta(x,z)$ depending only on $xz$ and satisfying
\beq
L_\beta (x, \partial_x) \Psi_\beta(x,z) =   z^N \Psi_\beta(x,z),
\label{2.17'2}
\eeq
where
\beq
L_{\beta}(x,\partial_x) = x^{-N} (D_x - \beta_1)
(D_x - \beta_2)\cdots(D_x - \beta_N),
\label{2.16}
\eeq
which is called a {\it{Bessel operator}} ($D_x = x \partial_x$). The
corresponding plane
$V_\beta \in\Gr$ in Sato's Grassmannian is called a {\em{Bessel plane\/}}
(it has an admissible basis
$w_{-k}=\partial^k_x\Psi_\beta(x,z)|_{x=1}$ if we take $v_k = e^z z^{-k}$ as a
basis of $\Vset$).
\ede
\bre{be}
In the above definition we use a convention from \cite{BHY3} which
we shall recall.
For a plane $W\in\Gr$ such that $\Psi_W(x,z)$ is well defined for
$x=x_0$  we set $v_k=e^{x_0 z} z^k$
and consider the subspace $W^{x_0}$ of $\Vset$ with an admissible basis
$w_{k}=\partial^k_x\Psi_W(x,z)|_{x=x_0}$.
The wave functions of $W^{x_0}$ and $W$ are connected by
$
\Psi_{W^{x_0}}(x,z) = e^{-x_0 z} \Psi_{W} (x + x_0, z)
$
and obviously
$$
\tau_{W} (t_1,t_2,t_3,\ldots) =
\tau_{W^{x_0}} (t_1-x_0,t_2,t_3,\ldots)
$$
where the RHS is considered as a formal power series in
$t_1-x_0,t_2,t_3,\ldots$

Throughout the paper we work with the spaces $W^1$
without explicitly mentioning it
and our tau-functions are formal power series in $t_1-1,t_2,t_3,\ldots$
\qed
\ere
Because the Bessel wave function depends only on $xz$, \eqref{2.17'2} implies
\beqa
&&D_x \Psi_\beta(x,z)= D_z \Psi_\beta(x,z),
\label{2.17'1}\\
&&L_\beta (z, \partial_z) \Psi_\beta(x,z) =
   x^N \Psi_\beta(x,z).
\label{2.17'3}
\eeqa

The monomial Darboux transformations of Bessel wave functions were introduced
in our previous paper \cite{BHY3}. They are a part of the solutions to
the bispectral problem (polynomial Darboux transformations) which we
constructed there.

First recall the definition of polynomial Darboux transformations given in
\cite{BHY3}.
\bde{dt}
We say that a plane $W$ (or the corresponding wave
function $\Psi_W(x,z)$) is a {\em{Darboux transformation\/}} of the
Bessel plane $V_\beta$ (respectively wave function $\Psi_\beta(x,z)$) iff there
exist polynomials $f(z)$, $g(z)$ and differential operators $P(x,\partial_x)$,
$Q(x,\partial_x)$ such that
\beqa
&&\Psi_W(x,z)=\frac{1}{g(z)} P(x,\partial_x) \Psi_\beta(x,z),
\label{3.8} \\
&&\Psi_\beta(x,z)=\frac{1}{f(z)} Q(x,\partial_x) \Psi_W(x,z).
\label{3.9}
\eeqa
The Darboux transformation is called {\em{polynomial\/}} iff
the operator $P(x, \partial_x)$ from \eqref{3.8} has the form
\beq
P(x,\partial_x)=x^{-n}\sum_{k=0}^n p_k(x^N) D_x^k,
\label{5.22}
\eeq
where $p_k$ are rational functions, $p_n\equiv 1$.

(In \cite{BHY3} we normalized $g(z)$.
For the present paper the normalization is unnecessary.)
\ede

There are two equivalent definitions of monomial Darboux transformations of
Bessel wave functions (see \cite{BHY3}).
\bde{mon1}
We say that the wave function $\Psi_W(x,z)$ (or the corresponding plane $W$)
is a {\em monomial Darboux transformation\/} of the Bessel wave function
$\Psi_\beta(x,z)$ (respectively the plane $V_\beta$) iff it is a polynomial
Darboux
transformation of $\Psi_\beta(x,z)$ with $g(z) f(z) = z^{d N}$, $d\in\Nset$.
\ede
\bde{mon2}
The wave function $\Psi_W(x,z)$ (or the corresponding plane $W$) is a
{\em monomial Darboux transformation\/} of the Bessel wave function
$\Psi_\beta(x,z)$ (respectively the plane $V_\beta$) iff \eqref{3.8} holds
with
$g(z)=z^n$, $n= \ord P$ and the kernel of the operator $P(x, \partial_x)$  has
a basis consisting of several groups of the form
\beq
\partial_y^l \Big(\sum_{k=0}^{k_0}\sum_{j=0}^{\mult(\beta_i+kN)-1}
b_{kj} x^{\beta_i + kN} y^j \Big) \Big|_{y = \ln x}, \quad 0 \leq l \leq j_0,
\label{mon}
\eeq
where $\mult(\beta_i + kN):=$ multiplicity of
$\beta_i+ kN$ in $\bigcup_{j=1}^N \{\beta_j+ N\Zset_{\geq0}\}$ and
$j_0 =\max\{j | b_{kj}\not=0 {\textrm{ for some }} k\}$.
\ede
Denote the set of monomial Darboux transformations of $V_\beta$ by
$\GrMBb$. For polynomial Darboux transformations we use the notation $\GrBb$.

\smallskip
A simple consequence of the above definitions is that
\beq
Q(x,\partial_x) P(x,\partial_x) = L_\beta(x,\partial_x)^d,
\label{c1}
\eeq
where
\beq
\ord P = n, \; \ord Q = dN-n \quad (g(z)=z^n, \; f(z) = z^{dN-n}).
\label{c2}
\eeq
Note that the monomial Darboux transformations have the following transitivity
and reflexivity properties:
\beqa
&&W\in\GrMBb, \; V_\beta\in\GrMB{\beta'}  \;\Rightarrow\; W\in\GrMB{\beta'};
\nn\\
&&V_\beta\in\GrMB{\beta'} \;\Leftrightarrow\; V_{\beta'}\in\GrMBb.
\nn
\eeqa
\subsection{}
In this subsection we recall the definition of $\W$, its subalgebras
$\WN$ and their bosonic representations introduced in \cite{BHY1}.

The algebra $w_{\infty}$ of the additional symmetries of the KP--hierarchy is
isomorphic to the Lie algebra of regular polynomial differential operators on
the circle
$$\D = \span \{ z^{\a} \partial_z^{\beta}| \; \a,  \beta
\in \Zset, \; \beta \geq 0 \}. $$
Its unique central extension \cite{KP1, KR} will be denoted by $W_{1+\infty}$.
This algebra gives the action of the additional symmetries on tau-functions
(see \cite{ASvM}). Denote by $c$ the central element of $W_{1+\infty}$ and by
$W(A)$  the image of $A \in \D$ under the natural embedding $\D \hookrightarrow
W_{1+\infty}$ (as vector spaces). The algebra $W_{1+\infty}$ has a basis
$$c, \; J_k^l = W(-z^{l+k} \partial_z^l), \qquad l,k \in \Zset, \; l \geq 0.$$

The commutation relations of $W_{1+\infty}$ can be written most conveniently in
terms of generating series \cite {KR}
\beq
\Bigl[ W(z^k e^{xD_z}),W(z^m e^{yD_z}) \Bigr] = ( e^{xm} - e^{yk}) W(z^{k+m}
e^{(x+y)D_z})+ \delta_{k,-m} \frac{e^{xm}-e^{yk}}{1-e^{x+y}}c,
\label{2.1}
\eeq
where $D_z = z\partial_z$.

Instead of working with the generators $J_k^l$ it is much more convenient to
work with the generating functions or fields (of dimension $l+1$)
\beq
 J^l(z) = \sum_{k \in \Zset } J_k^l z^{-k-l-1}. \label{2.3}
\eeq
The modes $J_k=J^0_k$ of the $\hat u(1)$ current $J(z) = J^0(z)$
generate the Heisenberg algebra:
\beq
\Bigl[ J_n, J_m \Bigr] = n\delta_{n,-m}c.
\label{heis}
\eeq
Recall its canonical representation in the bosonic Fock space $B$:
\beq
J_n= \frac{\partial}{\partial t_n},\quad
J_{-n}=nt_n,\quad n>0, \quad
J_0= Q\frac{\partial}{\partial Q}, \quad
c=1.
\label{heis1}
\eeq
It is well known that for $c=1$ the fields $J^l(z)$ can be expressed as
normally ordered polynomials in the current $J(z)$:
\beq
 J^l(z) =  l! \,\normprod{ S_{l+1} \Bigl( \frac{ J(z) }{ 1! },
 \frac{ \partial J(z) }{ 2! }, \ldots \Bigr) }.
\label{2.11}
\eeq
Here as usual
\[ \normprod{ J_n J_m } =
           \left\{\begin{array}{ll}
                    J_n J_m & \mbox{
                           $\textrm{for}\quad m>n$}\\
                    J_m J_n & \mbox{
                           $\textrm{for}\quad m<n$}
            \end{array}\right. \]
and the elementary Schur polynomials $S_l(t)$ are
determined by the generating series
\beq
\exp \Bigl( \sum_{k=1}^{\infty} t_k z^k \Bigr) =
   \sum_{l=0}^{\infty} S_l(t) z^l.
\label{2.10'8}
\eeq
Substituting \eqref{heis1} in \eqref{2.11} we obtain a bosonic representation
of
$\W$ with central charge $c=1$. For an explanation and a proof of \eqref{2.11}
see the Appendix.

In \cite{BHY1} we constructed a family of highest weight modules of $\W$ using
the above bosonic representation. We shall sum up the results from that paper
in a suitable for our purposes form. First we introduce the subalgebra $\WN$ of
$\W$ spanned by $c$ and $J_{kN}^l,$ $l,k\in \Zset,$ $l \geq 0 $. It is a simple
fact that $\WN$ is isomorphic to $\W$ (see \cite{FKN}).
\bth{2.5}
The functions $\tau_\beta(t)$ satisfy the constraints
    \beqa
&& J_0^l \tau_\beta = \lambda_\beta (J_0^l) \tau_\beta, \quad l\geq 0,
\label{2.19} \hfill\\
&& J_{kN}^l \tau_\beta = 0, \quad k>0, \; l\geq 0,
\label{2.20} \hfill\\
&& W\Bigl( z^{-kN} P_{\beta,k}(D_z) D_z^l \Bigr) \tau_\beta = 0,
\quad k>0, \; l\geq 0,
\label{2.21} \hfill
    \eeqa
where
$P_{\beta,k}(D_z) = P_\beta(D_z) P_\beta(D_z-N) \cdots P_\beta(D_z-N(k-1))$
and $P_\beta(D_z) = (D_z-\beta_1)\cdots(D_z-\beta_N)$.
\eth
The first two constraints mean that $\tau_\beta$ is a highest weight vector
with highest weight $\lambda_\beta$ of a representation of $\WN$ in the module
    \beq
\M_\beta = \span \Bigl\{ J_{k_1 N}^{l_1} \cdots J_{k_p N}^{l_p} \tau_\beta
\Big| k_1 \leq\ldots \leq k_p < 0 \Bigr\}.
\label{2.22}
    \eeq
In \cite{BHY1} we studied $\M_\beta$ as modules of $\W$. We proved that they
are
{\em{quasifinite\/}} (see \cite{KR}) and we derived formulae for the highest
weights and for the singular vectors. The latter formula turns
out to be the simplest corollary of \thref{7.4} (see \exref{7.6}).
\subsection{}
In the next sections we shall need the action of the so-called {\em{adjoint
involution\/}} $a$ on the modules $\M_\beta$.
On the tau-functions it acts as follows \cite{W}:
\beq
\tau_{aV}(t_1, t_2,\ldots,t_k,\ldots) =
\tau_V(t_1,-t_2,\ldots,(-1)^{k-1}t_k,\ldots).
\label{4.2}
\eeq
We continue this action on
$B^{(0)}=\Cset\left[\left[t_1-1, t_2, t_3,\ldots\right]\right]$
(cf.\ \reref{be}).

We shall continue it also on the elements of $\W$ in its bosonic representation
\eqref{2.11} naturally demanding
$$
a(U\tau) = a(U) a(\tau) \; {\textrm{  for }} U\in\W, \; \tau\in B^{(0)}.
$$
It acts on the Heisenberg algebra by $a(J_k) = (-1)^{k-1} J_k$, i.e.\
$a(J(z)) = J(-z)$, and on the fields $J^l(z)$ via \eqref{2.11}.
\bpr{ad}
If $\tau\in\M_\beta$ {\rm(}$\beta\in\Cset^N${\rm)}  then $a(\tau)\in\M_{a(\beta)}$ where
$a(\beta) = (N-1)\delta - \beta$, $\delta = (1,1,\ldots,1)$.
\epr
\proof
Using the commutation relations \eqref{2.1} one can prove by induction on the
dimension of the fields $J^l(z)$ that $a$ preserves $\WN$ (as a basis of the
induction one uses \eqref{2.11} for $l = 0, 1, 2$).  In the Appendix we give
another proof of this fact providing an explicit expression for a basis of $\W$
in which $a$ acts diagonally. The proposition now follows from the fact that
$a(\tau_\beta) = \tau_{a(\beta)}$ (see \cite{BHY3}).
\qed
\sectionnew{Tau-functions in Bessel modules as monomial Darboux
            transformations}
This section examines the connection between the class of representations
obtained in \cite{BHY1} (see Subsect.~1.3) and a part of the solutions to the
bispectral problem constructed in \cite{BHY3} (see Subsect.~1.2).
Our main result is the following.
\bth{7.1}
If $\tau_W$ is a tau-function lying in the $\WN$-module
$\M_\beta$ {\rm(}$\beta \in \Cset^N${\rm)} then the corresponding plane $W \in
\GrMBb$. Conversely, if $W \in \GrMBb$ then $\tau_W \in \M_{\beta'}$
for some $\beta' \in \Cset^N$ such that $V_{\beta'} \in \GrMBb$.
\eth
In general $\beta' \not= \beta$. A more precise version of the second part of
the theorem is given in Theorems \ref{t7.4} and \ref{t7.4a} below.
\subsection{}
For the {\em{proof\/}} of the first part of \thref{7.1} we shall need two
lemmas.
\ble{7.2}
If $\tau\in \M_\beta$ then
$\tau = u\cdot\tau_\beta$ with $u$ of the form
\beq
u=\sum a^l_k J^{l_1}_{-Nk_1}\cdots J^{l_r}_{-Nk_r},
\label{7.1}
\eeq
such that   all $l_i<Nk_i$.
\ele
\noindent
\proof
For $w=W\bigl(z^kP(D_z)\bigr)$ set $\rho(w)=\ord P+k$. Because of \thref{2.5}
for each  $w\in W_{1+\infty}(N)$ there exists $\widetilde w\in W_{1+\infty}(N)$
such that $w{\tau_\beta}=\widetilde w{\tau_\beta}$ and $\rho(\widetilde
w)<0$. Then for $w_1,\ldots, w_r\in W_{1+\infty}(N)$ we prove by induction on
$r$ that $w_1\cdots w_r\tau_\beta$    is a sum of     elements of the form
$\widetilde w_1\cdots \widetilde w_s \tau_\beta$ with $\rho(\widetilde w_i)<0$,
$s\le r$. Indeed, for $w\in W_{1+\infty}(N)$
\beqa
w\widetilde w_1\cdots \widetilde w_s \tau_\beta &=&
\widetilde w_1\cdots \widetilde w_s w\tau_\beta
+\left[w,\widetilde w_1\cdots \widetilde w_s\right] \tau_\beta
\nn\\
&=&
\widetilde w_1\cdots \widetilde w_s \widetilde w\tau_\beta +
\displaystyle\sum_{i=1}^s \widetilde w_1\cdots\left[w,\widetilde
w_i\right]\cdots \widetilde w_s \tau_\beta
\nn
\eeqa
with $\rho(\widetilde w)<0$.
\qed
\ble{7.3}
Let $X(t,z)$ be the vertex operator \eqref{1.8}.
Then
$$
X(t,z)J_k^l=\left(J_k^l +lJ_k^{l-1} +\delta_{l,0}\delta_{k,0} -z^{k+l}
\partial_z^l\right) X(t,z).
$$
\ele
\noindent
The {\it proof} of \leref{7.3} is given in the Appendix.

\medskip
\noindent
Now we can give the {\it proof} of the first part of \thref{7.1}. Let
$\tau_W=u\tau_\beta$ be a tau-function and $u$ be an element of the universal
enveloping algebra of $\WN$ of the form \eqref{7.1}. We compute the wave
function
$$ \Psi_W(x,z)=\frac{X(t,z)\tau_W(t)}{\tau_W(t)}\Big|_{t_1=x,\;
t_2=t_3=\cdots=0}. $$
Using \leref{7.3} we commute $X(t,z)$ and $u$ to obtain
$$
\Psi_W(x,z)=\frac{U(t,z)X(t,z)\tau_\beta(t)}{u\tau_\beta(t)}\Big|_{t_1=x,\;
t_2=t_3=\cdots=0}.
$$
where
$$
U(t,z)=\sum a^l_k
\left(J^{l_1}_{-Nk_1} + l_1J^{l_1-1}_{-Nk_1} -z^{-Nk_1+l_1} \partial_z^{l_1}
\right)\cdots
$$
$$
\cdots
\left(J^{l_r}_{-Nk_r} + l_rJ^{l_r-1}_{-Nk_r} -z^{-Nk_r+l_r} \partial_z^{l_r}
\right).
$$
{}From the bosonic formula \eqref{2.11} and the gradation of $W_{1+\infty}(N)$
it is clear that
$$
J^{l}_{-Nk}|_{t_1=x,\; t_2=t_3=\cdots=0}
=x^{Nk}\delta_{l+1, Nk}\quad {\rm if}\ \ l<Nk.
$$
What is relevant for us is that there are no differentiations in $t_1,
t_2,\ldots$ but only
a multiplication by powers of $x^N$. This gives for $U(t,z)$ the representation
$$
U(t,z)|_{t_1=x,\; t_2=t_3=\cdots=0}=
\sum a^l_k \left( x^{Nk_1}(\delta_{l_1+1, Nk_1} +l_1\delta_{l_1, Nk_1}) -
z^{-Nk_1+l_1}\partial_z^{l_1}\right) \cdots
$$
$$
\cdots \left( x^{Nk_r}(\delta_{l_r+1, Nk_r} +l_r\delta_{l_r, Nk_r}) -
z^{-Nk_r+l_r}\partial_z^{l_r}\right)
= z^{-mN} P(x^N,z^N,D_z),
$$
for some $m\in\Nset$ and a polynomial $P$ in $x^N$, $z^N$ and $D_z$.

In the same way $u|_{t_1=x,\; t_2=\cdots=0}=g(x^N)$ is polynomial in $x^N$.
Therefore
\beq
\Psi_W(x,z)=\frac{P(x^N,z^N, D_z)\Psi_\beta(x,z)}{z^{mN}g(x^N)}.
\label{7.2}
\eeq
Using (\ref{2.17'2}, \ref{2.17'1}) we obtain
\beq
\Psi_W(x,z) = z^{-mN} P_1(x^N,D_x)\Psi_\beta(x,z)
\label{7.3}
\eeq
for some operator $P_1$ with rational coefficients.

We also need an expression for $\Psi_\beta$ in terms of $\Psi_W$. It can be
obtained by using the adjoint involution $a$. By \prref{ad}
$\tau_{aW}=a(u)\tau_{a(\beta)}$
is a tau-function lying in the module $\M_{a(\beta)}$ and \eqref{7.3} gives
$$
\Psi_{aW}(x,z) = z^{-kN}P_2(x^N,D_x)\Psi_{a(\beta)}(x,z)
$$
for some operator $P_2$.
To complete the proof that $W\in\GrMBb$ we shall apply the following simple
lemma (see e.g.\ \cite{BHY3}, Proposition~1.7~(i)).
\ble{}
If the wave functions $\Psi_W(x,z)$ and $\Psi_V(x,z)$ satisfy
$$
\Psi_W(x,z)=\frac{1}{g(z)} P(x,\partial_x) \Psi_V(x,z),
$$
then
\beq
\Psi_{aV}(x,z)=\frac{1}{\check g(z)} P^*(x,\partial_x) \Psi_{aW}(x,z)
\label{4.7}
\eeq
where $\check g(z)=g(-z)$ and ``*'' is the formal adjoint {\rm(}i.e.\ the
antiautomorphism such that $\partial_x^* = - \partial_x$, $x^*= x${\rm)}.
\ele
The above lemma leads to
\beq
\Psi_\beta(x,z) = z^{-kN}(-1)^{kN} P_2^*(x^N,D_x) \Psi_W(x,z)
\label{7.3b}
\eeq
which combined with \eqref{7.3} proves that $W \in \GrMBb$.
\subsection{}
The second part of \thref{7.1} is a consequence of Theorems \ref{t7.4} and
\ref{t7.4a} below. Before stating the first of them let us introduce some
notation and recall some simple facts.
\ble{5.1}
Let $\beta \in \Cset^N$ and $\alpha \in \Cset^M$. Then

{\em(i)} $L_\alpha L_\beta=L_\gamma$,  where
$$
\gamma=(\alpha_1+N,\alpha_2+N,\ldots,\alpha_M+N,
\beta_1,\beta_2,\ldots,\beta_N);
$$

{\em(ii)} $(L_\beta)^d=L_{\beta^d}$,  where
$$
\beta^d=(\beta_1,\beta_1+N,\ldots,\beta_1+(d-1)N,\ldots,\beta_N,\ldots,\beta_N
+(d-1)N);
$$

{\em(iii)}  If
$\{\beta_1,\ldots,\beta_N\}=\{\underbrace{\alpha_1,\ldots,\alpha_1}_{k_1}
,\ldots,
\underbrace{\alpha_s,\ldots,\alpha_s}_{k_s} \}$
with distinct $\alpha_1, \ldots,\alpha_s$, then
$$
\ker L_\beta=\span\left\{ x^{\alpha_i}(\ln x)^k\right\}_{1\le i\le s,\ 0\le
k\le k_i-1}.
$$
\ele
\noindent
The {\em{proof\/}} is obvious.

\smallskip

Let $W\in Gr_{MB}(\beta)$ be a monomial Darboux transformation of the Bessel
plane $V_\beta$, $\beta\in \Cset^N$. We can consider only the case when
$n\le d$ (see eqs.\ (\ref{c1}, \ref{c2})) since the general case can be
reduced to this one by  a
left multiplication of $Q$ by $L_\beta$, which does not change $W$ and
$\beta$. Let $\gamma = \beta^d$ (see \leref{5.1} (ii)), i.e.\
\beq
\gamma_{(k-1)d+j} := \beta_k + (j-1)N, \quad 1\le k\le N,\; 1\le j\le d.
\label{7.4'}
\eeq
First we consider the case when $\ker P$ has a basis of the form
\beq
f_k(x) = \sum_{i=1}^{dN} a_{ki} x^{\gamma_i},\quad
0\le k\le n-1,\quad
\label{7.4}
\eeq
i.e.\ there are no logarithms. \deref{mon2} in this case is equivalent to
\beq
\gamma_i-\gamma_j \in N\Zset \setminus 0
\quad {\rm if}\ a_{ki}a_{kj}\not=0, i\not=j.
\label{7.4'5}
\eeq
We say that the element $f_k(x)$ of the above basis of $\ker P$ is
{\em{associated}} to $\beta_s$ ($1\leq s \leq N$) iff
\beq
\gamma_i-\beta_s\in N\Zset_{\ge0}\; {\rm if}\; a_{ki}\not=0.
\label{7.7}
\eeq
Then up to a relabeling we can take a subset $\{ \beta_s\}_{1 \leq s\leq M}$
such that
\beq
\beta_s-\beta_t\not\in N\Zset\quad {\rm for}\ \  1\le s\not=t \le M
\label{7.6}
\eeq
and each element of the basis \eqref{7.4} of $\ker P$ is associated to some
$\beta_s$ from this set. Denote by $n_s$ the number of elements associated
to $\beta_s$ and set $n_s=0$ for $s>M$.
Then $n_1+\cdots+n_N = n$. We put
\beq
\beta'=(\beta_1+n_1N-n, \beta_2+n_2N-n,\ldots,\beta_N+n_NN-n).
\label{7.8}
\eeq
\bth{7.4}
Let $W$ be a monomial Darboux transformation of the Bessel plane $V_\beta$
with $\ker P$ satisfying {\rm(\ref{7.4}, \ref{7.4'5})} and $\beta'$ be as
above.
Then the tau-function  $\tau_W$ of $W$ lies in the $W_{1+\infty}(N)$-module
$\M_{\beta'}$.
\eth
\noindent
{\it Proof}.
We shall use the following formula
\beqa
\Psi_W(x,z)&=&\frac{Wr\bigl( f_0(x),\ldots, f_{n-1}(x),\Psi_\beta(x,z)\bigr)}
{z^n Wr\bigl( f_0(x),\ldots, f_{n-1}(x)\bigr)}
\nn\\
&=&\frac{\sum\det A^I Wr\bigl(x^{\gamma_I}\bigr) \Psi_I(x,z)}
{\sum\det A^I Wr\bigl(x^{\gamma_I}\bigr)},
\hfill \label{3.26}
\eeqa
where $Wr$ denotes the Wronski determinant.
The sum is  taken over all $n$-element subsets
$I=\{i_0<i_1<\ldots <i_{n-1}\} \subset\{0,1,\ldots, dN-1\}$,
$x^{\gamma_I} = \{x^{\gamma_i}\}_{i\in I}$,
$A^I=(a_{k,i_l})_{0\le k,\; l\le n-1}$
and $\Psi_I(x,z)$ is the
wave function of the above type of Darboux transformations with
$f_k(x) = x^{\gamma_{i_k}}$, i.e.\
$$
\Psi_I(x,z) = z^{-n} L_{\gamma_I}(x,\partial_x) \Psi_\beta(x,z),
\qquad \gamma_I = (\gamma_{i_0},\ldots,\gamma_{i_{n-1}}).
$$
It is important also that
\ble{5.2}{\rm{\cite{BHY3}}}
$\Psi_I(x,z)$  is again a Bessel wave function{\rm:}
\beq
\Psi_I(x,z)=\Psi_{\gamma+dN\delta_I -n\delta} (x,z),
\label{5.6}
\eeq
where the vectors $\delta_I$, $\delta$ are defined by
$$
(\delta_I)_i=\cases{1, &if $i\in I$\cr
                    0, &if $i\not\in I$\cr}
$$
 and
$$
\delta_i=1\ \ {\rm for\ all}\ \ i\in\{1,\ldots, dN\}.
$$
\ele
We set
\beq
I_0 := \{1,\ldots, n_1, d+1,\ldots, d+n_2, \ldots, (N-1)d+1,\ldots,
(N-1)d+n_N\}.
\label{7.9}
\eeq
Then (see \eqref{7.4'})
$$
\gamma_{I_0}=\{ \beta_1,\beta_1+N,\ldots, \beta_1+(n_1-1)N, \ldots, \beta_N,
\beta_N+N,\ldots, \beta_N+(n_N-1)N\}
$$
and clearly
\beq
\tau_{I_0}=\tau_{\beta'}
\label{7.11}
\eeq
(recall that $\tau_\beta=\tau_\gamma$ when $L_\beta^d=L_\gamma$).

First we shall consider the case when $\det A^{I_0}\not= 0$.
Without loss of generality we can put $\det A^{I_0}=1$.
Let $A_0$ be the $n\times dN$ matrix
$\bigl(a_{ki}\delta_{i,i^0_k}\bigr)_{0\le k\le n-1,\; 1\le i\le dN}$ where
$I_0=\{i_0^0<i_1^0<\ldots< i^0_{n-1}\}$ is from \eqref{7.9}.
For $\zeta\in \Cset$ we define the matrix $A(\zeta)$ as follows
\beq
A(\zeta)=\zeta A+(1-\zeta)A_0.
\label{7.12}
\eeq
Then $A(\zeta)_{ki} =a_{ki}$ for $i=i_k^0$ and ${}=\zeta a_{ki}$ for $i\not=
i_k^0$. Thus \eqref{7.4'5} holds with $A(\zeta)_{ki}$ instead of $a_{ki}$ and
the Darboux
transformation $W(\zeta)$ of $V_\beta$ with a matrix $A(\zeta)$ is monomial:
$$
W(\zeta)\in Gr_{MB}(\beta).
$$
The main idea of the proof of \thref{7.4} is to consider $W(\zeta)$ as a
deformation of $W(0)=V_{\beta'}$. We shall prove that $\tau_{W(\zeta)}\in
\M_{\beta'}$ for all $\zeta$, hence $\tau_W=\tau_{W(1)}\in \M_{\beta'}$.
We first need a lemma expressing $\Psi_I$ in terms of $\Psi_{\beta'}\equiv
\Psi_{I_0}$.
\ble{7.5}
If $\det A(\zeta)^I\not=0$ for some $\zeta$ then
\beq
\Psi_I(x,z)=x^{-q_I}Q_I(z,\partial_z)\Psi_{I_0}(x,z),
\label{7.13}
\eeq
where $Q_I$ is a Bessel operator of order $q_I$, divisible by $N$ and
satisfying
\beq
q_I\le p_I := \sum_{i\in I} \gamma_i- \sum_{i\in I_0} \gamma_i.
\label{7.14}
\eeq
The number $p_I$ is also divisible by $N$.
\ele
\noindent
{\it Proof}. For $1\le s\le M$ we set $I_s=\{i\in I\mid \gamma_i-\beta_s\in
N\Zset_{\ge0}\}$. Then (\ref{7.7}, \ref{7.6}) imply that
$$
I=\bigcup_{s=1}^M I_s,\quad I_s\cap I_t=\emptyset\quad
{\textrm{ for }}\;  s\not= t \;{\textrm{  and }}\;
\card I_s=n_s.
$$
Let
$$
\gamma_{I_s}=\left\{ \beta_s+\nu_1^{(s)}N, \beta_s+(\nu_2^{(s)}+1)N,
\ldots, \beta_s+(\nu_{n_s}^{(s)}+n_s-1)N\right\}
$$
   for $1\le s\le M$ and
$\gamma_{I_s}=\emptyset$ for $s>M$,
where $\nu_i^{(s)}\in \Zset$,
$0=\nu_1^{(s)}=\cdots=\nu_{k_s}^{(s)}< \nu_{k_s+1}^{(s)}\le \ldots \le
\nu_{n_s}^{(s)}$. Then
$$
p_I=\sum_{s=1}^M \sum_{i=1}^{n_s} N\nu_i^{(s)}\ge
N\sum_{s=1}^M(n_s-k_s)=N(n-k),
$$
where $k=\sum_{s=1}^M k_s$; we set $k_s=0$ for $s>M$.
We take $q_I=N(n-k)$ and $Q_I=L_\alpha$ be the Bessel operator of order $q_I$
such that
\beq
L_\alpha L_{\gamma_{I_0}} =L_{\gamma_I}(L_\beta)^{n-k}.
\label{7.15}
\eeq
We shall prove that such $L_\alpha$ exists. Using \leref{5.1} the right hand
side of \eqref{7.15} can
be represented as
$$
L_{\gamma_I}(L_\beta)^{n-k}=L_{\alpha'},
$$
where
\beqa
\alpha'&=&\bigl(\gamma_I+(n-k)N \delta_I \bigr) \cup \beta^{n-k}
\nn\\
&=& \bigcup\limits_{s=1}^N\Bigl\{ \bigl(\gamma_{I_s}+(n-k)N \delta_{I_s}\bigr)
\cup \bigl(\beta_s,\beta_s+N,\ldots,\beta_s+(n-k-1)N\bigr)\Bigr\}.
\nn
\eeqa
We see that $\alpha'$ includes
$$
\beta_s,\beta_s+N,\ldots, \beta_s+(n-k+k_s-1)N,\quad (1\le s\le N)
$$
and therefore includes $\gamma_{I_0}$, which proves \eqref{7.15}. Now the proof
of \eqref{7.13} is straightforward. Using that
$$
\Psi_I=x^{-n} L_{\gamma_I}\Psi_\beta,\quad
\Psi_{I_0}=x^{-n} L_{\gamma_{I_0}}\Psi_\beta
$$
(the Bessel operators act in the variable $z$), we compute
\beqa
x^{-q_I} Q_I\Psi_{I_0} &=& x^{-(n-k)N} L_\alpha x^{-n}
L_{\gamma_{I_0}}\Psi_\beta =x^{-(n-k)N-n} L_\gamma(L_\beta)^{n-k}\Psi_\beta
\nn\\
&=& x^{-(n-k)N-n} L_\gamma x^{(n-k)N}\Psi_\beta
=x^{-n}L_\gamma\Psi_\beta=\Psi_I.
\nn
\eeqa
\qed

\smallskip

Now we can apply formula \eqref{3.26}. Obviously
\beq
Wr\bigl(x^{\gamma_I}\bigr)=\Delta_I x^{\sum_{i\in I}\gamma_i -\frac{n(n-1)}2},
\label{7.16}
\eeq
where for $I=\{ i_0<\ldots<i_{n-1}\}$ we set
\beq
\Delta_I=\prod_{r<s}(\gamma_{i_r}-\gamma_{i_s}).
\label{7.17}
\eeq
Using this and \eqref{7.13} we write $\Psi_{W(\zeta)}$ as
$$
\Psi_{W(\zeta)}(x,z)=\frac{\sum \det A(\zeta)^I \Delta_I x^{p_I-q_I}
          Q_I(z,\partial_z) \Psi_{I_0}(x,z)}
{\sum \det A(\zeta)^I \Delta_I x^{p_I}}.
$$
We expand the denominator around $\zeta=0$.
Using that $\Psi_{I_0}=\Psi_{\beta'}$, that $p_I$ and $q_I$ are divisible by
$N$ and \eqref{2.17'3} for $\beta'$  we obtain
\beq
\Psi_{W(\zeta)}(x,z)=\sum_{i\ge 0} \zeta^i
P_i(z,\partial_z)\Psi_{\beta'}(x,z)
\label{7.19}
\eeq
for some operators $P_i$ without constant term for $i\geq 1$ and with $P_0
\equiv 1$. Indeed \eqref{7.12} implies
that $A(\zeta)^{I_0}=A^{I_0}$ and $\det A(\zeta)^{I_0}=1$ does not depend on
$\zeta$. Therefore for $i\ge 1$ $P_i$ is a linear combination of operators
$$
Q_I(L_{\beta'})^{(kp_I-q_I)/N},\quad k\geq 1, \; I\not= I_0
$$
which are nontrivial Bessel operators (see \leref{7.5}) and thus do not have a
constant term. For $I=I_0$ $Q_{I_0}\equiv1$ and $p_I=q_I$, now $\det
A(\zeta)^{I_0}=1$ implies $P_0\equiv 1$.
Denoting $w_{-k}=\partial_x^k\Psi_{\beta'}(x,z)|_{x=1}$ we obtain (see
(\ref{tau}, \ref{1.2a'}))
$$
\begin{array}{rl}
\tau_{W(\zeta)}
&=\sigma\Bigl\{\left(\sum \zeta^iP_iw_{0}\right)\wedge
\left(\sum\zeta^iP_iw_{-1}\right)\wedge\ldots\Bigr\}
\\
\noalign{\vskip3pt}
&=\sigma\Bigl\{w_{0}\wedge w_{-1}\wedge w_{-2}\wedge\ldots +\zeta(
P_1w_{0}\wedge w_{-1}\wedge w_{-2}\wedge\ldots
\\
\noalign{\vskip3pt}
&\qquad{}+w_{0}\wedge P_1w_{-1}\wedge w_{-2}\wedge\ldots+\cdots)+\cdots\Bigr\}
\\
\noalign{\vskip3pt}
&=\tau_{\beta'} +\zeta r(P_1)\tau_{\beta'} +
\zeta^2\bigl(r(P_2)+\frac12 r(P_1)^2 -\frac12
r(P_1^2)\bigr)\tau_{\beta'}+\cdots
\end{array}
$$
We see that all coefficients at the powers of $\zeta$ are polynomials in
$r(P_i^k)$ applied to $\tau_{\beta'}$ and thus belong
to the $W_{1+\infty}(N)$-module $\M_{\beta'}$ (see \eqref{2.2'5}). Now we shall
use the formula
\beq
\tau_{W(\zeta)}=
\frac{\sum \det A(\zeta)^I \Delta_I \tau_I }
{\sum \det A(\zeta)^I \Delta_I}
\label{3.28}
\eeq
(see \cite{BHY2}). Because the numerator depends polynomially on $\zeta$ the
above considerations show that it belongs to $\M_{\beta'}$.
Setting $\zeta=1$ we obtain $\tau_W\in\M_{\beta'}$.
This completes the proof of \thref{7.4} in the case when $\det A^{I_0}\not=0$.

The general case can be deduced again from the fact that the numerator of
\eqref{3.28} is
 polynomial in the entries of $A$. Up to a relabeling one can suppose that the
first $n_1$ functions of the basis \eqref{7.4} of $\ker P$ are associated to
$\beta_1$, the next $n_2$ to $\beta_2$, etc. Then the Darboux transformation
with a matrix
$$
A(\xi) = A + \xi E_0, {\textrm{where }}
E_0 = (\delta_{i\,i^0_k})_{1\le i\le dN,\,0\le k\le n-1}
$$
is monomial (see \eqref{7.4'5}). Obviously $\det A^{I_0} = \det (A^{I_0} +
\xi E) \not= 0$ for all but a finite number of $\xi\in\Cset$
(where $E$ is the identity matrix) and for
them $\tau_{W(\xi)} \in\M_{\beta'}$. Because the numerator of \eqref{3.28}
(with $\zeta$ replaced with $\xi$) is a polynomial in $\xi$, it belongs to
$\M_{\beta'}$ for all $\xi\in\Cset$ and for $\xi=0$ it is exactly $\tau_W$
(recall that a tau-function is defined up to a multiplication by a constant).
\qed
\subsection{}
Now we shall consider the general case of a monomial Darboux transformation
of $V_\beta$, $\beta \in \Cset^N$. Using repeatedly \leref{5.2} with $n=d=1$ we
see that $V_\beta \in \GrMB{\nu}$ with $\nu$ of the form
\beq
\nu= (\underbrace{\nu_1,\nu_1, \ldots,\nu_1}_{N_1},\ldots,
      \underbrace{\nu_p,\nu_p, \ldots,\nu_p}_{N_p})
\label{7.100}
\eeq
such that
\beq
\nu_i-\nu_j \not\in N \Zset \quad {\rm{for}} \; i \not= j
\label{7.101}
\eeq
$(N_1+\cdots+N_p=N).$

Let $W\in \GrMB{\nu}$, i.e.\ $ \ker P$ has a basis consisting of several groups
of the form described in \deref{mon2}:
\beq
\sum_{j=l}^{j_0} l! {j \choose l} f_j(x)
(\ln x)^{j-l}, \quad 0 \leq l \leq j_0
\label{7.102}
\eeq
where $j_0 \leq N_i-1$ and
\beq
f_j(x) = \sum_{k = 0}^{d-1} b_{kj} x^{\nu_i + kN}.
\label{7.102'}
\eeq
We say that the element
\eqref{7.102} of $\ker P$ has {\em{level}} $j_0 - l$. For $1 \leq s \leq p$ we
denote by $n_s^r$ the number of elements in the basis of $\ker P$ of level $r$
associated to $\nu_s$ (see \eqref{7.7}). Put
\beq
\nu' = (\nu_1+n_1^0 N-n,\ldots,
        \nu_1+n_1^{N_1-1} N-n,\ldots,
        \nu_p+n_p^0 N-n,\ldots,
        \nu_p+n_p^{N_p-1} N-n).
\label{7.103}
\eeq
\bth{7.4a}
If $W$ is a monomial Darboux transformation of $V_\nu$ with $\nu$, $\nu'$ and
$\ker P$ as above then $\tau_W \in \M_{\nu'}$.
\eth
\proof We shall make a limit in \thref{7.4}. We note that since
$(\ln x)^j = \partial_\lambda^j x^\lambda |_{\lambda = 0}$ we have
$$(\ln x)^j = \lim_{\epsilon \to 0} \epsilon^{-j} \sum_{k=0}^{j}
(-1)^{k} {j \choose k} x^{-\epsilon k}.$$
Set
$\nu(\epsilon) = (\nu_1,\nu_1+\epsilon,\ldots,\nu_1+(N_1-1)\epsilon,
                 \ldots, \nu_p,\nu_p+\epsilon,\ldots,\nu_p+(N_p-1)\epsilon)$.
Consider the Darboux transformation $W(\epsilon)$ of $V_{\nu(\epsilon)}$ with a
basis of $\ker P(\epsilon)$ consisting of groups of the form
(cf.\ (\ref{7.102}, \ref{7.102'}))
$$
g_l(x)= x^{\epsilon (j_0-l)} \sum_{j=l}^{j_0} l! {j \choose l} f_j(x)
\sum_{k=0}^{j-l} \epsilon^{l-j}
(-1)^k {j-l \choose k} x^{-\epsilon k}, \quad 0 \leq l \leq j_0.
$$
We shall show that this transformation is monomial. More precisely,
we shall prove that $\ker P(\epsilon)$ has a basis consisting of groups of
elements of the form
$$
h_l(x)= x^{\epsilon (j_0-l)}
\sum_{j=l}^{j_0} l! {j \choose l} \epsilon^{l-j} f_j(x),
\quad 0 \leq l \leq j_0.
$$
This is an obvious consequence of the identity
$$g_l(x) = \sum_{k=0}^{j_0-l} \frac{(-1)^k \epsilon^{-k}}{k!} h_{l+k}(x),
\quad 0 \leq l \leq j_0.$$
We apply \thref{7.4} for $W(\epsilon)$ noting that exactly $n_s^r$ elements of
the above basis of $\ker P(\epsilon)$ are associated to $\nu_s + r\epsilon$.
Taking the limit $ \epsilon \to 0$ completes the proof.
\qed
\subsection{}
As an illustration to \thref{7.4} we shall consider the case $n=d=1$.
Now the matrix $A$ is $1\times N$:
$$
A=(a_1\ a_2\ \ldots\ a_N),
$$
subsets $I$ consist of one element: $I=\{i\}$, and
\beq
\Psi_I\equiv\Psi_i=\frac1z\left(\partial_x-\frac{\beta_i}x\right)\Psi_\beta
=\Psi_{(\beta_1-1,\ldots,\beta_i+N-1,\ldots,\beta_N-1)}.
\label{7.20}
\eeq
The formula \eqref{3.26} now becomes
\beq
\Psi_W(x,z)=\frac1{zx}\left(x\partial_x-\frac{\sum a_i\beta_ix^{\beta_i}}{\sum
a_i x^{\beta_i}}\right) \Psi_\beta(x,z).
\label{7.21}
\eeq
Let $a_1 \not= 0$. The Darboux transformation is monomial when
$$
\beta_i-\beta_1=N\alpha_i,\quad \alpha_i\in\Zset\quad {\rm for}\ \ a_i\not=0
$$
and up to a relabeling we can suppose that all $\alpha_i$ are positive.
Then
$
I_0=\{1\}$, $\Psi_{I_0}\equiv\Psi_1\equiv\Psi_{\beta'}$, where
$\beta'=(\beta_1+N-1,\beta_2-1,\ldots,\beta_N-1)$.
Using that
\beqa
\Psi_\beta(x,z) &=& x^{-N} L_\beta(z,\partial_z)\Psi_\beta(x,z) = (xz)^{-N}
P_\beta(D_z)\Psi_\beta(x,z)
\nn\\
&=&
{\displaystyle{(xz)^{-N+1} \frac{P_\beta(D_z+1)}{D_z+1-\beta_1}
\Psi_{\beta'}(x,z) }}
\nn
\eeqa
and that
\beq
P_\beta(D_z)=\frac{D_z-\beta_1}{D_z-(\beta_1+N)} P_{\beta'}(D_z-1),\quad
P_{\beta'}(D_z)=\frac{D_z-(\beta_1+N-1)}{D_z-(\beta_1-1)} P_\beta(D_z+1)
\label{7.22}
\eeq
we obtain from \eqref{7.21}
\beq
\Psi_W(x,z)=\frac1{\sum a_ix^{N\alpha_i}} P_1(z,\partial_z) \Psi_{\beta'}(x,z),
\label{7.23}
\eeq
where
\beq
P_1(z,\partial_z) =\sum_{i=1}^N a_i
\left\{ -N\alpha_i z^{-N}\frac{P_{\beta'}(D_z)} {D_z-\beta_1'}\left(z^{-N}
P_{\beta'} (D_z)\right)^{\alpha_i}\right\}.
\label{7.24}
\eeq
Then $A(\zeta)=(a_1\ \zeta a_2\ \ldots\ \zeta a_N)$, the numerator of
\eqref{3.28} is equal to $\tau_{\beta'}+\zeta r(P_1) \tau_{\beta'}$ and up to a
constant
\beq
\tau_W = r(P_1) \tau_{\beta'}.
\label{7.26}
\eeq
\bex{7.6}
Let $A=(0\ 1\ 0\ \ldots\ 0)$ and
$\beta'_2-\beta_1'=N\alpha=N(\alpha_2-1)$, $\alpha\in\Zset_{\ge0}$.
Set
$$
\beta''=(\beta_1'-N, \beta_2'+N, \beta_3', \ldots, \beta_N')
=(\beta_1-1,\beta_2+N-1,\beta_3-1,\ldots,\beta_N-1).
$$
Then the module $\M_{\beta''}$ embeds in $\M_{\beta'}$. The singular
vector $\tau_{\beta''}$ is given by (\ref{7.26}, \ref{7.24}).
Up to some changes of notation ($\beta_i=Nr_i$, etc.) in this way we recover
Theorem 6 from \cite{BHY1}.
\qed
\eex
\bex{even}
Let $N=2$, $\beta_1 -\beta_2 =2\a$, $\a \in\Zset_{\ge0}$. Then the tau-function
$\tau_\a:=\tau_{(1/2 +\a, 1/2-\a)}$ is highest weight vector for the reducible
$\W(2)$-module
$\M_{(1/2+\a, 1/2-\a)}$, which will be denoted below by $\M_\a$.
\exref{7.6} gives that $\M_0\supset\M_2\supset\M_4\supset\ldots$ and
$\M_1\supset\M_3\supset\M_5\supset\ldots$. Any bispectral tau-function
corresponding to an ``even'' potential \cite{DG} can be obtained by a monomial
Darboux transformation with $d=n\le \a$ from $\tau_{\a}$
as shown  in \cite{MZ} (see also \cite{BHY3}, Example 5.3). \thref{7.4} shows
that $\tau_W\in\M_{\a-n}$. On the other hand
$\tau_W\in\Gr^{(2)}$, which gives that $\tau_W$ belongs to the $Vir$-module
$M_{\a-n}^\infty$  introduced in \cite{HH} (see also \cite{F}). The
modules $M_{\alpha}^\infty$, $\alpha\in\Zset_{\ge0}$ are shown to be the
reducible
Verma modules over $Vir$ with $c=1$ (whose highest weight vectors are
the above tau-functions $\tau_\alpha$). In this way we obtain:
\eex
\bco{even2}
Any tau-function $\tau_W$ of an ``even'' potential can be obtained by a
monomial
Darboux transformation from the highest weight vector $\tau_\a$ of a reducible
$Vir$-module $M_\alpha^\infty$, $\alpha\in\Zset_{\ge0}$ {\rm(}defined in
{\rm\cite{HH})} and it belongs
to
the module $M_{\alpha-n}^\infty$, where $n$ is the order of the
operator $P$ {\rm(}$n\le\a${\rm)}.
Conversely, any tau-function in  $M_\alpha^\infty$ is tau-function of an
``even'' potential {\rm\cite{DG}}.
\eco

Consider the set of modules $\M_\beta$ of the most degenerate case
$\beta_i - \beta_j \in N\Zset$ for all $i,j$ ($\beta \in \Cset^N$). The
embeddings among these modules are described by $N$ lattices: the $k$-th one of
them having a maximal module $\M_{\beta^{(k)}}$,
$$\beta^{(k)} = (\underbrace{b+N, \ldots, b+N}_{k},
                 \underbrace{b, \ldots, b}_{N-k}),\;
                 b=\frac{N-2k-1}{2},$$
for $0 \leq k \leq N-1$ (cf.\ \exref{7.6}). \leref{5.2} implies that the set of
monomial Darboux transformations of $\beta^{(k)}$ with $\ord P\in N \Zset$
coincides with the set of monomial Darboux transformations of
$\beta^{(0)}$ with $\ord P\in k+N \Zset$. In the latter case the corresponding
$\M_{\nu'}$ given by \thref{7.4a} belongs to the $k$-th lattice, i.e.\ it is a
submodule of $\M_{\beta^{(k)}}$.
So we obtain the following corollary.
\bco{11}
The manifold of monomial Darboux transformations from $\beta^{(k)}$ with $\ord
P \in N \Zset$ coincides with the manifold of tau-functions lying in the module
$\M_{\beta^{(k)}}$.
\eco
\bre{}
We shall consider another aspect of \thref{7.4}. Recall that the subalgebras
$\W(d)$, $d\in\Nset$ of $\W$ are isomorphic to $\W\equiv\W(1)$ and a
representation of $\W(d)$ with central charge $N$ gives rise to a
representation of $\W$ with central charge $dN$. Each singular vector of
$\W$ is obviously a singular vector of $\W(d)$ but the converse is not
always true. It is an intersting question to describe the latter.

In our terminology this question can be reformulated as follows.
{\em{Which Bessel
tau-functions $\tau_\alpha$, $\alpha\in\Cset^{dN}$ lie in a $\WN$-module
$\M_\beta$, for some $\beta\in\Cset^N$}}?
Such tau-functions are given by
\thref{7.4}: if $I$ is an $n$-element subset of $\{1,\ldots,dN\}$
then $\tau_{\beta^d + dN\delta_I - n\delta} \in\M_{\beta'}$ (see \eqref{7.8}
and \leref{5.2}).

Let us consider the simplest case $N=1$ and set $d=2n$,
$I=\{1,3,\ldots,2n-1\}$. Then $\beta=(0)$, $\beta^d = (0,1,\ldots,2n-1)$ and
$$
\beta^d + dN\delta_I - n\delta = (1-n,3-n,\ldots,n-1) \cup (n,n+2,\ldots,3n-2)
= (1-n,n)^n.
$$
So we obtain that $\tau_{(1-n,n)}$ lies in $\W(1)$ module
$\M_{(0)} = \Cset[t_1,t_2,\ldots]$, i.e.\ they are polynomials of
$t_1,t_2,\ldots$ (obviously it coincides with the module over the Heisenberg
algebra with highest weight vector $\tau_{(0)}=1$). These
tau-functions play an important role in the ``KdV case''
of \cite{DG} and are connected with the rank $1$ bispectral algebras
which contain an operator of order $2$. The general case of rank $N$
bispectral algebras containing an operator of order $dN$ motivates the study
of the above considered ``embeddings''.
\ere
\sectionnew{Hierarchies of symmetries of the manifolds of monomial Darboux
transformations}
An immediate consequence of the results of the previous section is the
existence
of hierarchies of symmetries preserving the manifolds of monomial Darboux
transformations.

Denote by $\W^+(N)$ the subalgebra of $\WN$ spanned by $J_{Nk}^l$, $k,l \geq
0$.
\bth{7.7}
For $\beta\in\Cset^N$ the vector fields corresponding to $W_{1+\infty}^+(N)$
are tangent to the manifold $\GrMBb$ of monomial Darboux transformations.
More precisely, if $W\in \GrMBb$ then
\beq
\exp\left(\sum_{i=1}^p \lambda_i J_{Nk_i}^{l_i}\right)\tau_W
\label{7.27}
\eeq
is a tau-function associated to a plane from $\GrMBb$ for
arbitrary
$p\in \Nset$, $\lambda_i\in\Cset$, $l_i,k_i\in\Zset_{\ge0}$.
\eth
\noindent
\proof
Indeed, if
$W\in Gr_{MB}(\beta)$, $\beta\in \Cset^N$ then $\tau_W\in \M_{\beta'}$ for some
$\beta'\in \Cset^N$ such that $V_{\beta'} \in \GrMBb$. Now from the gradation
of $W_{1+\infty}(N)$ it is clear that
$\sum_{i=1}^p \lambda_i J_{Nk_i}^{l_i}$ acts nilpotently on $\tau_W$,
i.e.\ \eqref{7.27} is well-defined and also belongs to the module
$\M_{\beta'}$.
Moreover, it is a tau-function (see \cite{ASvM}) and \thref{7.1} shows that its
corresponding plane belongs to $Gr_{MB}(\beta')= \GrMBb$.
\qed

\smallskip

Let us introduce the following terminology.
A $\beta \in \Cset^N$ is called generic if $V_\beta$ cannot be obtained by a
Darboux transformation of some $V_\alpha$, $\alpha \in \Cset^M$, $M<N$.
We put $\GRMBN = \bigcup_\beta \GrMBb$,
$\beta\in\Cset^N$--generic. These manifolds
are important because they give bispectral algebras of rank $N$ \cite{BHY3}.
Therefore \thref{7.7} implies:
\bco{12}
The manifold $\GRMBN$ of monomial Darboux transformations which give bispectral
algebras of rank $N$ is preserved by the vector fields corresponding to
$\W^+(N)$.
\eco
An interesting question is when a polynomial Darboux transformation of a Bessel
operator $L_\beta$ of order $N$ gives again an operator of order $N$ (see
\cite{DG} for $N=2$ and \cite{BHY3}). In \cite{BHY3}, Proposition 5.4
we proved that for generic
$\beta$ such transformation is necessarilly monomial.
The corresponding manifold $\GrMBb\cap Gr^{(N)}$ is
also preserved by an hierarchy of symmetries.
More precisely, in the bosonic realization \eqref{2.11} of $W_{1+\infty}(N)$ we
put $J_{kN}=0$, $k\in \Zset$ and define
$$
\overline L_m=\frac1N J_{mN}^1|_{J_{kN}=0,\;k\in\Zset}
=\frac1N\sum_{i\in{\Zset\setminus N\Zset}} \normprod{J_{mN-i}J_i}.
$$
The operators $\overline L_m$, $m\in \Zset$ form a Virasoro algebra with
central charge $N-1$  which we denote by $Vir_N$.
Denote by $Vir_N^+$ the subalgebra spanned by $\overline L_m$, $m\ge0$. Then we
can formulate the following theorem which for $N=2$ contains Magri--Zubelli's
result \cite{MZ}.
\bth{7.8}
The manifold $\GrMBb\cap Gr^{(N)}$ is preserved by
the vector fields corresponding to $Vir_N^+$. More precisely, if $W\in
\GrMBb\cap Gr^{(N)}$ then
\beq
\exp\left(\sum_{i=1}^p \lambda_i\overline L_{k_i}\right) \tau_W
\label{7.28}
\eeq
is a tau-function associated to a plane from $\GrMBb\cap Gr^{(N)}$ for
arbitrary $p\in \Nset$, $\lambda_i\in\Cset$, $k_i\ge0$.
\eth
\noindent
{\it Proof}. Obviously formula \eqref{7.28} gives exactly the same result as
$$
\exp\left(\sum \lambda_i N^{-1} J_{k_i N}^1\right) \tau_W.
$$
This is because $\tau_W(t)$ does not depend on $t_{kN}$ and in $J_{k_i N}^1$,
$k_i\ge0$ the variables $t_{kN}$ are present only as coefficients of
differentiations with respect to $t_{mN}$. This implies that \eqref{7.28} is a
tau-function of a plane belonging to $Gr^{(N)}$.
\qed

\bre{3.2,5}
The manifold $\GrMBb\cap\Gr^{(dN)}$, $\beta\in\Cset^N$ is
preserved
by the vector fields corresponding to the subalgebra of $\W^+(N)$ generated by
\beqa
&& J_{kN}, \quad {\textrm{for }} d \not| k, \, k\ge 0;
\nn\\
&& {\widetilde L}_m = \sum_{i\in\Zset\setminus\ dN\Zset}
\normprod{J_{mN-i} J_i},
\quad {\textrm{for }} m\ge 0.
\nn
\eeqa
(Because on $\Gr^{(dN)}$ ${\widetilde L}_m$ acts as $J^1_{mN}$.)

For $N=1$, $d=2$ we recover the well known fact that the potentials from the
``KdV  case'' of \cite{DG} are preserved by the KdV flows.

For generic $\beta\in\Cset^N$ $\GrMBb\cap\Gr^{(dN)}$ gives bispectral
algebras of rank $N$ containing an
operator of order $dN$. However for $d>1,$ $N>1$ these are not all such
algebras, cf.\ \cite{BHY3}. It is still an open problem to describe
the symmetries of the latter.
\qed
\ere

We shall conclude this section with some comments. As $Gr^{(N)}$ is a reduction
of $Gr$, the (associative) algebra $W_N$ is a reduction of $W_{1+\infty}(N)$
-- see e.g.\ \cite{FKN, vM}. In more details, the fields $J^0(z),
J^1(z),\ldots,J^{N-1}(z)$ generate the (vertex operator) algebra ${\mathcal
W}(gl_N)$ \cite{FKRW}. Its reduction ${\mathcal W}(sl_N)$ is obtained by
putting
$J_{kN}=0$, $k\in\Zset$ in \eqref{2.11}; the modes of the corresponding fields
generate the so-called $W_N$ algebra. (More precisely, this is
a representation of $W_N$ with $c=N-1$.) Then $Vir_N\subset W_N$ and we
conjecture that \thref{7.8} is valid with $W_N^+$ instead of $Vir_N^+$.
\section*{Appendix}
\cleqn
\renewcommand{\theequation}{{\rm{A}}.\arabic{equation}}
In this appendix we give the technical proofs of some of the results from
Sect.~1 and explain them in more details.

First, following \cite{KRa}, we shall recall the boson-fermion correspondence.
Recall the definition of the fermionic Fock space $F$ from subsection 1.1.
The free fermions can be realized as wedging and contracting operators:
\begin{eqnarray*}
&&\psi_{-j+\frac12}\left( v_{i_0}\wedge v_{i_1}\wedge\ldots\right)
= v_j\wedge v_{i_0}\wedge v_{i_1}\wedge \ldots\\
&&\psi^*_{j-\frac12}\left( v_j\wedge v_{i_0}\wedge v_{i_1}\ldots\right)
= v_{i_0}\wedge v_{i_1}\wedge\ldots
\end{eqnarray*}
They satisfy the canonical anticommutation relations
\beq
\left[\psi_\lambda, \psi^*_\mu\right]_+ =\delta_{\lambda,-\mu},
\quad
\left[\psi_\lambda, \psi_\mu\right]_+ =0,
\quad
\left[\psi^*_\lambda, \psi^*_\mu\right]_+ =0,
\label{1.1}
\eeq
where $[a,b]_+ = ab+ba$.

Let $gl_\infty$ be the Lie algebra of all $\Zset\times \Zset$ matrices having
only a finite
number of non-zero entries. One can define a representation $r$ of $gl_\infty$
in the fermionic Fock space $F$ as follows. For the basis $E_{ij}\in gl_\infty$
put
\beq
r\left(E_{ij}\right)=\psi_{-i+\frac12}\psi^*_{j-\frac12}
\label{1.2a}
\eeq
and continue this by linearity. Then for $A\in gl_\infty$
\beq
r(A) (w_0\wedge w_{-1}\wedge w_{-2}\wedge\ldots) =
Aw_0\wedge w_{-1}\wedge w_{-2}\wedge\ldots +
w_0\wedge Aw_{-1}\wedge w_{-2}\wedge\ldots + \cdots \, .
\label{1.2a'}
\eeq
The above defined representation $r$ obviously cannot be continued on
the Lie algebra $\widetilde{gl}_{\infty}$ of all $\Zset\times \Zset$ matrices
with finite number of non-zero diagonals. If we regularize it by
\beq
\hat r\left(E_{ij}\right)=
\normprod{\psi_{-i+\frac12}\psi^*_{j-\frac12}}, \label{1.2}
\eeq
where as usual $\normprod{\psi_\mu\psi_\nu^*}=\psi_\mu\psi^*_\nu$ for
$\nu>0$
and $-\psi^*_\nu\psi_\mu$ for $\nu<0$, this will give a representation for the
central extension $\widehat{gl}_\infty =\widetilde{gl}_\infty \oplus \Cset c$
of
$\widetilde{gl}_\infty$. Here the central charge $c$ acts as a multiplication
by 1. Define the free fermionic fields
$$
\psi(z)=\sum_{j\in\Zset} \psi_{j-\frac12} z^{-j}
\quad{\rm and}\quad
\psi^*(z)=\sum_{j\in\Zset} \psi^*_{j-\frac12} z^{-j}.
$$
Then the anticommutation relation \eqref{1.1} can be written as
\beq
\Bigl[ \psi(z_1),\psi^*(z_2) \Bigr]_+=\delta(z_{12}),
\label{1.3}
\eeq
where $z_{12}=z_1-z_2$ and $\delta(z_{12})=\sum_{n\in\Zset}
z^n_1z^{-n-1}_2$.
Introduce also the $\hat u(1)$ current
\beq
J(z)=\normprod{\psi^*(z)\psi(z)}=\sum_{n\in \Zset} J_n z^{-n-1}.
\label{1.4}
\eeq
The modes $J_n$ generate the Heisenberg algebra \eqref{heis}.

The above introduced spaces $F^{(m)}$ are spaces of irreducible representations
of the Heisenberg algebra with charge $m$ and central charge $c=1$.
Using that such a representation is unique up to isomorphism we obtain the
isomorphism known as the boson-fermion correspondence (\ref{bo-fe},
\ref{heis1}).
In terms of the states $|m\rangle =v_m\wedge v_{m-1}\wedge\ldots$
and the operator $H(t)=-\sum_{k=0}^\infty t_kJ_k$ we have for $|\varphi \rangle
\in F$
\beq
\sigma(\,|\varphi\rangle) =\sum_{m\in\Zset} \langle
m|\;e^{H(t)}\;|\varphi \rangle Q^m.
\label{1.5}
\eeq
We also introduce the scalar bosonic field:
\beq
 \phi(z) = \hat q + J_0 \log z+ \sum_{n \not=0} J_n \frac{z^{-n}}{-n}
 \label{2.5}
\eeq
with operator product expansion
$ \phi(z_1) \phi(z_2) \sim \log (z_1-z_2),$
which is equivalent to \eqref{heis} and
\beq
\Bigl[ J_n , \hat{q} \Bigr] = \delta_{n,0},
\label{2.5'5}
\eeq
 and such that
$$
\exp \hat q = Q, \quad
J(z) = \partial_z \phi (z), \quad
Q^m = \normprod{ e^{m \phi(z)}} |0\rangle|_{z=0}.
$$
Then the fermionic fields $\psi (z)$, $\psi^{*} (z)$ act on the bosonic Fock
space $B$ as
\beq
 \psi^{*} (z) = \normprod{e^{\phi (z)}}, \qquad
 \psi (z) = \normprod{e^{-\phi (z)}}.
 \label{2.6}
\eeq
Here as usual $\normprod{J_n J_m} = J_n J_m$
for $m>n$, $\normprod{J_n J_m} = J_m J_n$ for $m<n$ and
$\normprod{\hat q J_0} = \normprod{J_0 \hat q} = \hat q J_0 $.

One can define a natural embedding of $W_{1+\infty}$ in $\widehat{gl}_{\infty}$
in the following way \cite{KP1, KR}. Consider a realization of $\Vset$ as
the space of Laurent series in $z^{-1}$. Fixing the basis $v_j = z^{-j}$ of
$\Vset$ each element $A \in \D$ corresponds to a matrix
$\phi_0 (A) \in \widetilde{gl}_{\infty}$, i.e.\ one defines an embedding
$ \phi_0\colon \D \hookrightarrow \widetilde{gl}_{\infty}$.
This embedding can be extended to an embedding
$$ \hat \phi_0\colon W_{1+\infty} \hookrightarrow  \widehat{gl}_{\infty}.$$
In the case when $c$ acts as a multiplication by $1$ we can obtain, using
\eqref{1.2}, a free field realization of $W(A)$:
\beq
 W(A) = \res_{z=0} \normprod{ \psi(z) A \psi^{*}(z)}
 \label{2.2}
\eeq
for $ A \in \D$. In the notation of \eqref{1.2} this means
\beq
W(A) = \hat r(A).
\label{2.2'5}
\eeq
Also note that $ \hat r(A) = r(A) $ for operators $A$ having $ \diag
\phi_0(A) = 0$.

{}From \eqref{2.2} we derive a bosonic realization of the fields $J^l(z)$:
\beq
 J^l(z) = \normprod{ \bigl( \partial_z^l \psi^{*}(z) \bigr) \psi(z) }
 \label{2.7}
\eeq
which combined with \eqref{2.6} gives
\beq
 J^l(z) = \frac{1}{l+1} \normprod{ e^{-\phi(z)} \partial_z^{l+1} e^{\phi(z)} }.
 \label{2.8}
\eeq
This and the Taylor formula imply \eqref{2.11}. Indeed,
\beqa
  && \sum_{l \geq 0} \frac{x^{l+1}}{l!} J^l(z) =
      \normprod{ \Bigl( \sum_{l \geq 0} \frac{x^{l+1}}{(l+1)!}
      \partial_z^{l+1} e^{\phi(z)} \Bigr)
      e^{-\phi(z)} } \nn \\
  &&= \normprod{ \Bigl( e^{\phi(z+x)}-e^{\phi(z)} \Bigr) e^{-\phi(z)} }
      = \normprod{ \Bigl( e^{ \sum_{k \geq 0} \frac{x^k}{k!}
      \partial^k \phi(z)} e^{-\phi(z)}-1 \Bigr)} \nn \\
  && = \sum_{l \geq 1} x^l \normprod{ S_l \Bigl( \frac{ \partial \phi }{1!},
      \frac{ \partial^2 \phi }{ 2! }, \ldots \Bigr)}. \nn
\eeqa
Comparing the coefficients at $x^{l+1}$ and using that
$J(z) = \partial_z \phi(z)$ we get \eqref{2.11}.

To describe the action of the involution $a$ on $\W$ we introduce another
convenient basis $V^l_k,\, c$ of $\W$ through the fields
\beq
 V^l(z) = \sum_{k \in \Zset } V_k^l z^{-k-l-1}. \label{2.10'5}
\eeq
These fields are  {\em{quasiprimary\/}} of dimension $l+1$ with respect to the
Virasoro algebra
generated by $V^1_k$, $c$ and are used essentially by Cappelli, Trugenberger
and Zemba in their study of the Quantum Hall Effect (see \cite{CTZ} and
references therein). They can be defined by \cite{BGT}
\beq
 V^l(z) = \frac{l!}{(2l)!} \partial_1^l (-\partial_2)^l \Bigl\{ z_{12}^l
 \normprod{\psi^{*}(z_1)\psi(z_2)} \Bigr\}\Big|_{z_1=z_2=z}
\label{2.10'6}
\eeq
(for the connection with the fields $J^l(z)$ see \cite{BGT}, eq.\ (1.41)).
We have an analog of \eqref{2.8}:
\beq
 V^l(z) = \frac{1}{l} {2l \choose l}^{-1} \sum_{k=0}^{l-1}
{l \choose k}{l \choose k+1} \partial_1^{l-k} (-\partial_2)^{k+1}
\normprod{e^{\phi(z_1) - \phi(z_2)}}\Big|_{z_1=z_2=z}.
\label{2.10'7}
\eeq
This leads to an analog of \eqref{2.11} proved in the same way:
\beqa
V^l(z) &=& \frac{ (l-1)!\, l! \, (l+1)! }{ (2l)! }
\sum_{k=0}^{l-1} {l \choose k}{l \choose k+1}{l+1 \choose k+1}^{-1} (-1)^{k+1}
\nn\\
&\times& \normprod{
S_{l-k} \Bigl( \frac{J(z)}{1!}, \frac{\partial J(z)}{2!}, \ldots \Bigr)
S_{k+1} \Bigl( -\frac{J(z)}{1!}, -\frac{\partial J(z)}{2!}, \ldots \Bigr)
}  \label{2.12}
\eeqa
Substituting $a(J(z))=J(-z)$ for $J(z)$ in \eqref{2.12} it is easy to see that
\beq
a(V^l(z)) = V^l(-z),  \qquad l \ge 0
\label{av}
\eeq
(we use that the elementary Schur polynomials $S_l$ are homogeneous of degree
$l$ if $\deg t_k=k$).
In terms of the modes
$
a(V^l_k) = (-1)^{l+k+1} V^l_k
$
showing that $\WN$ is preserved by the involution $a$.

At the end we shall give the {\em{proof of \leref{7.3}}}.
Comparing (\ref{2.5},  \ref{2.6}) with \eqref{1.8} we see that
$$
\psi^*(z)=\normprod{ e^{\phi(z)} } = QX(t,z).
$$
{}From \eqref{2.7} and \eqref{1.3} we derive commutation relations
\beq
 \Bigl[ J^l(z_1),\psi^{*}(z_2) \Bigr] =  \lim_{z_3 \to z_1} \partial_{z_1}^l
 \normprod{ \Bigl[ \psi^{*}(z_1) \psi(z_3) ,\psi^{*}(z_2) \Bigr] }
                =  \delta (z_{12}) \partial^l_{z_2} \psi^{*}(z_2).
 \label{2.10}
\eeq
Using that $J_0 Q=Q(J_0+1)$ (see \eqref{2.5'5}) we compute
$$
Q^{-1}J^l(z) = Q^{-1} \lim_{z_{1,2}\to z} \frac{\partial_1^{l+1}}{l+1}
\normprod{ e^{\phi(z_1)-\phi(z_2)} }
=\lim_{z_{1,2}\to z} \frac{\partial_1^{l+1}}{l+1}
\normprod{ e^{\phi(z_1)-\phi(z_2)} }
\frac{z_1}{z_2} Q^{-1}
$$
$$
{}=\cases{(J^l(z)+lz^{-1}J^{l-1}(z))Q^{-1}, &for $l>0$.\cr
(J^0(z) +z^{-1})Q^{-1}, &for $l=0$.\cr}
$$
Then
$$
X(t,z_1) J^l(z_2) =Q^{-1}\psi^*(z_1) J^l(z_2)=Q^{-1} J^l(z_2) \psi^*(z_1) -
Q^{-1}\delta(z_{12})\partial_1^l\psi^*(z_1)
$$
$$
{}=\left( J^l(z_2) +z_2^{-1} l J^{l-1}(z_2)\right) X(t,z_1)
-\delta(z_{12})\partial_1^l X(t,z_1)
$$
(for $l=0$ instead of $lJ^{l-1}(z_2)$ put 1). Comparing the coefficients of
$z_2^{-k-l-1}$ in both sides completes the proof.
\begin{small}
    
\end{small}
\end{document}